\documentclass{article}
\usepackage[dvips,a4paper,portrait,verbose,margin=2cm]{geometry} 
\usepackage[cp1251]{inputenc} %еуееъшщ
\usepackage{oldgerm}
\usepackage{graphicx}
\usepackage{hyperref}

\usepackage[Verbose]{parallel}
\newcommand{\probb}[1]{{\bf #1}}

\newcommand{\numb}[1]{{\mbox{\sc #1}}}
\newcommand{\numbm}[1]{$\overline{\mbox{\sc #1}}$}

\newcommand{\threel}[3]{\begin{Parallel}{}{} \ParallelLText{#1} \ParallelRText{#2} \end{Parallel}}
\newcommand{\six}[6]{%
    \begin{Parallel}{}{}% 
         \ParallelLText{#1} \ParallelRText{#2}% 
         \ParallelPar%
         \ParallelLText{#4} \ParallelRText{#5}%
\end{Parallel}}

\newcounter{formula}

\newcommand{\fraccc}[2]{\leavevmode\kern.1em\raise.5ex\hbox{\scriptsize $#1$}
  \kern-.1em / \kern-.15em \lower.25ex\hbox{\scriptsize $#2$}\,}

\begin{document}

\title{Alcuin's Propositiones de Civitatibus:\\
                  the Earliest Packing Problems}
                          
\author{Nikolai Yu. Zolotykh}

%\date{August 5, 2013}
\date{Version 2: December 28, 2017}

\maketitle   

\begin{abstract}
We consider three problems about cities from Alcuin's {\it Propositiones ad acuendos juvenes}.
These problems can be considered as the earliest packing problems difficult also for modern state-of-the-art packing algorithms.
We discuss the Alcuin's solutions and give the known (to the author) best solutions to these problems.
\end{abstract}

\section{Introduction}

The manuscript {\it Propositiones ad acuendos juvenes}
({\it Problems to Sharpen the Young})
attributed to Alcuin of York is considered as the earliest collection of mathematical problems in Latin \cite{HadleySingmaster1992}.

In Alcuin's Opera Omnia published by Frobenius Forster \cite{Forster1777} and 
revised and republished by J.-P. Migne \cite{Migne1863} there is a version of {\it Propositiones} containing $53$ problems. 
Another version of {\it Propositiones} with three additional problems was published in Venerable Bede's Opera Omnia \cite{Herwagen1563};
it was revised and republished by J.-P. Migne \cite{Migne1862}.
The attribution of this version to Bede is considered to be spurious. 
A modern edition of {\it Propositiones} was made by Menso Folkerts \cite{Folkerts1978}. 
He found $12$ manuscripts, the oldest one is dating from the late $9$th century. It is incomplete, but already contains the three additional problems of Bede's text.

There are two translations of {\it Propositiones} into English. One of them was made by John Hadley with commentary by David Sigmaster and John Hadley \cite{HadleySingmaster1992}. A literal translation with commentary was made by Peter J. Burkholder \cite{Burkholder1993}.
Also, there is a German translation by Helmuth Gericke with commentary by Menso Folkerts and Helmuth Gericke \cite{FolkertsGericke1993}.
%A Critical German edition is 

According to Singmaster, {\it Propositiones} contain
``first River Crossing 
Problems ($3$ types); first Explorer's Problem; first Division of Casks; first 
Apple-sellers' Problem; first Collecting Stones; unusual solution of 
Posthu\-mous Twins Problem; first Three Odds Make an Even; first Strange 
Families'' \cite{Singmaster1996}.

On the other hand, some geometric problems from {\it Propositiones} with probably Roman origin use rather crude approximations
(for example, some of them require $\pi=4$).
Nevertherless certain geometric problems are of particular interest.

In {\it Propositiones} there are $11$ geometric problems,
IX, X, XXI--XXV, XXVII--XXI, dealing with the area of figures (triangles, rectangles, circles and quadrangles). 
Moreover, problems XXVII--XXVIII are formulated as packing problems. So they can be considered as the earliest packing problems.
They are all about cities (of different shapes) in which it is required to put rectangular houses of equal known sizes.
Problem XXVIII dealing with triangular city is easier than XXVII and XXVIIII about quadrangular and round cities, respectively.
And the formers are difficult even for modern state-of-the-art packing algorithms.
We give the original problem formulations and solutions of the problems with the translation by Burkholder \cite{Burkholder1993}
and give the best known solutions. 
%We don't discuss here the Alcuin's solutions because for these problems they all are incorrect (for discussion see \cite{HadleySingmaster1992,Burkholder1993}). 
%In discussion of the Alcuin's solutions we use also \cite{HadleySingmaster1992,Burkholder1993}.   

\section{Propositiones de civitatibus (problems about cities)}

\subsection{Quadrangular city problem}

\six{
\probb{XXVII. Propositio de civitate quadrangula.}
Est civitas  quadrangula quae habet in uno latere pedes mille centum; et in
alio latere pedes mille; et in fronte pedes \numb{dc}, et in altera pedes \numb{dc}. Volo
ibidem tecta domorum ponere, sic, ut habeat unaquaeque casa in longitudine
pedes \numb{xl},  et in latitudine pedes \numb{xxx}.  Dicat, qui velit, quot casas capere
debet?
}{
\probb{27. Proposition concerning the quadrangular city.}
There is  a quadrangular city which has one side of $1100$ feet, another side
of $1000$ feet, a front of $600$ feet, and a final side of $600$ feet.  I want to
put some  houses there so that each house is $40$ feet long and $30$ feet wide.
Let him say, he who wishes, How many houses ought the city to contain?
}{
%\probb{Задача о четырехугольном городе.}
%Есть четырехугольный город, который имеет одну сторону в $1100$ футов; противоположную --- $1000$ футов;
%во фронте --- $600$ футов и с противоположной стороны --- $600$ футов.
%Я хочу заложить в нем дома, так, что каждый дом $40$ футов в длину и $30$ в ширину.
%Пусть скажет, кто желает, сколько домов должен город вмещать.
}{
\probb{Solutio.}
Si fuerunt  duae  hujus  civitatis  longitudines  junctae,  facient  \numbm{ii}  \numb{c}.
Similiter duae,  si fuerunt  latitudines junctae,  faciunt \numbm{i}  \numb{cc}.  Ergo duc
mediam de  \numbm{i} \numb{cc},  faciunt [Bede:  fiunt] \numb{dc}, rursusque duc mediam de \numbm{ii} \numb{c}, fiunt \numbm{i} \numb{l}.
Et quia  unaquaeque domus  habet in  longitudine [Bede:  \dots{}in longo]  pedes \numb{xl}, et in lato
\numb{xxx}: deduc [Bede:  duc] quadragesimam partem de mille \numb{l}, fiunt \numb{xxvi}.  Atque iterum
assume tricesimam de \numb{dc}, fiunt \numb{xx}.  Vicies ergo \numb{xxvi} ducti, fiunt \numb{dxx}.  Tot
domus capiendae sunt.
}{
\probb{Solution.}
If the  two lengths  of this  city were joined together, they would measure
$2100$ [feet].   Likewise,  if the  two sides were joined, they would measure
$1200$.  Therefore, take half of $1200$, i.e. $600$, and half of $2100$, i.e. $1050$.
Because each  house is $40$ feet long and $30$ feet wide, take a fourtieth part
of $1050$,  making $26$.  Then, take a thirtieth of $600$, which is $20$.  $20$ times
$26$ is $520$, which is the number of houses to be contained in the city.
}{
%\probb{Решение.}
%Если сложить две стороны этого города, то получилось бы $2100$. 
%Аналогично, если сложить две ширины, то получилось бы $1200$.
%Следовательно, произведи половину от $1200$, получается $600$, 
%снова произведи половину от $2100$, получается $1050$.
%И так как каждый дом имеет длину $40$ футов и ширину $30$ футов, произведи сороковую часть от $1050$, получается $26$.
%И потом возьми тридцатую часть от $600$, получается $20$. 
%Итак, $26$ двадцатикратно, получается $520$.
%Столько домов вмещается.
}

\ 
\medskip

Alcuin uses the Egyptian (or Roman) formula for calculating the area of a quadrilateral 
$$
S=\frac{a+c}{2} \cdot \frac{b+d}{2},
$$
where $a$, $b$, $c$, $d$ are length of four sides of the quadrilateral
(the product of the half sum of its opposite sides).
According it, $S = 1050\times 600 = 630,000$  sq ft, that is equal to $525$ house areas.
If the shape of the city is an isosceles trapezium then its area is $627,808.689\dots$ sq ft, that is equal to $523.174\dots$ house areas.

Seemingly Alcuin is not interested in whether the city really could contain the number of houses indicated by him,
though discarding the fractional part in dividing $1050$ by $40$ could purpose to take into account this.

\begin{figure}
\centering
\includegraphics[height=7cm]{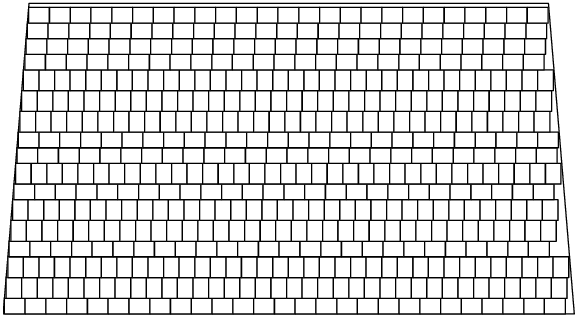} 
\caption{Problem XXVII. One can put $510$ houses in the trapezoidal city.}
\label{fig_quadrangula}
\end{figure}

Singmaster \cite{HadleySingmaster1992} notes that he can get $516$ houses in (if the shape of the city is an isosceles trapezium);
but ``by turning the houses the other way'' he can get $517$ in. 
``By having some houses either way'' he can get $519$ in.
I don't know these solutions and I can fit only $510$ houses in; see Fig.\,\ref{fig_quadrangula}.

\subsection{Triangular city problem}

\six{
\probb{XXVIII. Propositio de civitate triangula.}
Est civitas  triangula quae  in uno habet latere pedes \numb{c}, et in alio latere
pedes \numb{c},  et in  fronte  pedes  \numb{xc},  volo  enim  ibidem  aedificia  domorum
construere [Bede:  Volo ut fiat ibi domorum constuctio\dots{}], sic tamen, ut unaquaeque domus habeat in longitudine pedes
\numb{xx}, et  in latitudine  pedes \numb{x}.  Dicat, qui potest, quot domus capi debent?
}{
\probb{28. Proposition concerning the triangular city.}
There is  a triangular city which has one side of $100$ feet, another side of
$100$ feet,  and a  third of  $90$ feet.   Inside  of this,  I want  to build a
structure of  houses, however,  in such a way that each house is $20$ feet in
length, $10$  feet in width.  Let him say, he who can, How many houses should
be contained [within this structure]?
}{
%\probb{Задача о треугольном городе.}
%Есть треугольный город, который имеет по одной стороне $100$ футов, и по другой стороне $100$ футов, а по фронту --- $90$ футов,
%хочу же там расположить дома [Bede: Хочу создать структуру домов], так, чтобы каждый дом имел в длину $20$ Футов, а в ширину $10$ футов.
%Пусть скажет, кто может, сколько домов должно вмещаться [в этом городе].
}{
\probb{Solutio.}
Duo igitur  hujus civitatis latera juncta fiunt \numb{cc}, atque duc mediam de \numb{cc},
fiunt \numb{c}.   Sed quia in fronte habet pedes \numb{xc}, duc mediam de \numb{xc}, fiunt \numb{xlv}.
Et quia  longitudo uniuscujusque  domus habet pedes \numb{xx}, et latitudo ipsarum
pedes \numb{x}, duc \numb{xx} partem in [Bede:  de] \numb{c}, fiunt \numb{v}.  Et pars decima quadragenarii \numb{iv}
sunt.   Duc itaque  quinquies \numb{iiii},  fiunt \numb{xx}.  Tot domos hujusmodi captura [Bede: capienda] est civitas.
}{
\probb{Solution.}
Two sides  of the  city joined  together make $200$; taking half of $200$ makes
$100$.   But because  the front  is $90$ feet, take half of $90$, making $45$.  And
since the  length of  each house  is $20$ feet while the width is $10$, take $20$
into $100$,  making five.   A  tenth part of $40$ is four; thus, take four five
times, making $20$.  The city is to contain this many houses in this way.
}{
%\probb{Решение.}
%Две стороны этого города, сложенные вместе, дают $200$, и потом произведи половину от $200$, получается $100$.
%Но поскольку по фронту имеет $90$ футов, проиведи половину от $90$, получается $45$.
%И поскольку длина каждого дома имеет $20$ футов, а ширина $10$ футов, раздели $100$ на $20$ частей, получается $5$.
%И десятая часть от $40$ [sic!] есть $4$. Таким образом, возьми $4$ пять раз, получается $20$.
%Столько таких домов должен содержать город.
}

\medskip

\begin{figure}
\centering
\includegraphics[height=7.5cm]{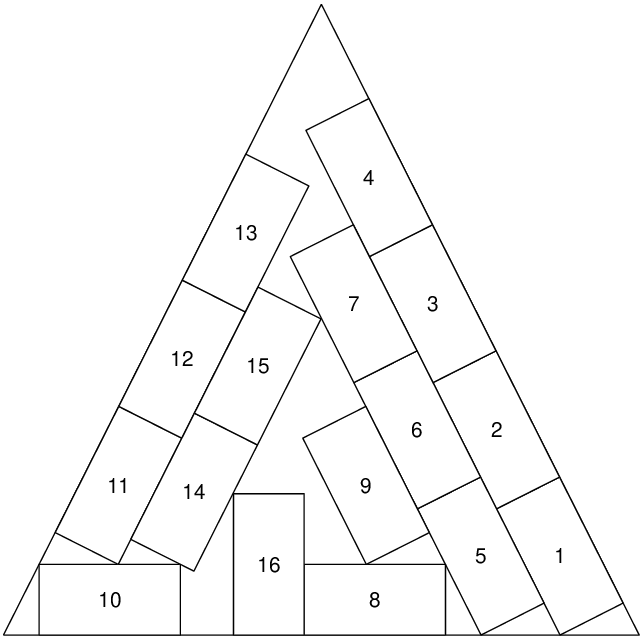} 
\caption{Problem XXVIII. One can put $16$ houses in the triangular city.}
\label{fig_triangula}
\end{figure}

\begin{figure}
\centering
\includegraphics[height=7.5cm]{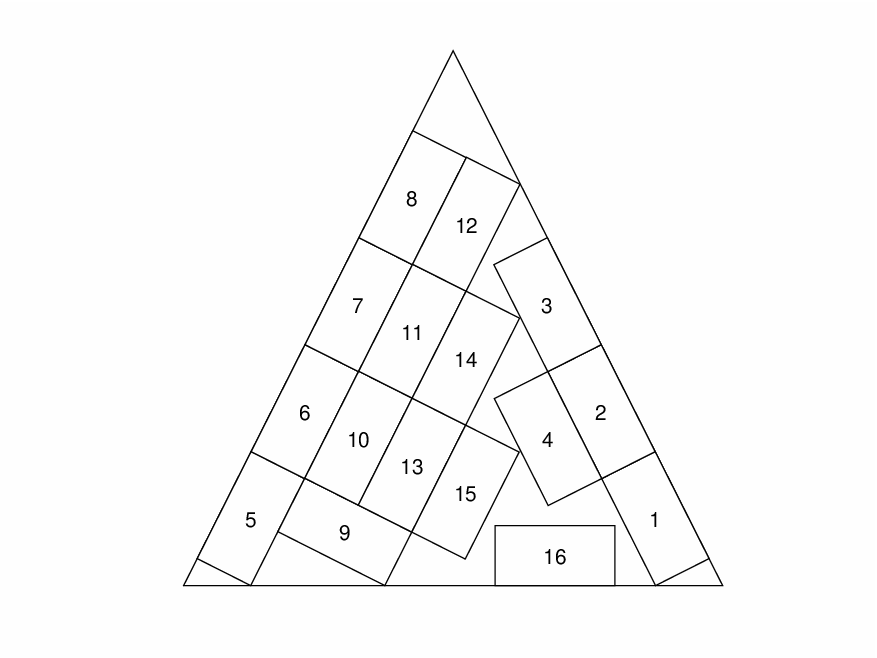} 
\caption{Problem XXVIII. Another way to put $16$ houses in the triangular city \cite{Hemme2010}.}
\label{fig_triangula2hemme}
\end{figure}

For calculation the area Alcuin uses the Egyptian formula that for triangles has the form
$$
S=\frac{a+b}{2} \cdot \frac{c}{2}.
$$
According it, the area of the city must be $100\times 45 = 4500$ sq ft, i.e. $22\frac{1}{2}$ house areas.
The correct value of the area is $4018.628\dots$, i.e. $20.093\dots$ house areas.
As Singmaster notes \cite{HadleySingmaster1992}, 
``The conversion of $45$ to $40$ may be an attempt to compensate for the inaccuracy of the formula, or to
allow for the difficulty of fitting rectangular houses into the triangular town.''
Singmaster can get $15$ houses in. Hadley can get $18$ houses in if ``the walls can be bent slightly'' \cite{HadleySingmaster1992}.
Here we give a solution with $16$ houses (and with straight walls). 
It is shown in Fig.\,\ref{fig_triangula}. 
Note that there is a small gap between houses $7$ and $15$.
Also, there is a gap between houses $8$ and $16$ or/and between houses $14$ and $16$.
Another solution to this problem with $16$ houses is proposed by H.\,Hemme \cite{Hemme2010};
see Fig.\,\ref{fig_triangula2hemme}.

\subsection{Round city problem}

\six{
\probb{XXVIIII. Propositio de civitate rotunda.}
Est civitas  rotunda quae  habet in circuitu pedum \numb{viii} millia.  Dicat, qui
potest, quot  domos capere  debet, ita  ut unaquaeque habeat in longitudine
pedes \numb{xxx}, et in latitudine pedes \numb{xx}?
}{
\probb{29. Proposition concerning the round city.}
There is  a city  which is $8000$ feet in circumference.  Let him say, he who
is able, How many houses should the city contain, such that each [house] is
$30$ feet long, and $20$ feet wide?
}{
%\probb{Задача о круглом городе.} 
%Есть круглый город, который имеет в окружности $8000$ футов.
%Пусть скажет, кто может, сколько домов он должен содержать, таких, что
%каждый имеет $30$ футов в длину и $20$ в ширину.
}{
\probb{Solutio.}
In hujus  civitatis ambitu  \numb{viii} millia  pedum numerantur, qui sesquialtera
proportione dividuntur  in \numbm{xxxx} \numb{dccc}, et  in  \numbm{iii} \numb{cc}.    In  illis  autem
longitudo domorum;  in istis latitudo versatur.  Subtrahe itaque de utraque
summa medietatem, et remanent de majori \numbm{ii} \numb{cccc}:  de minore vero \numbm{i} \numb{dc}.  Hos
igitur \numbm{i} \numb{dc} divide  in vicenos  et invenies  octoagies viginti,  rursumque
major  summa,   id  est   \numbm{ii} \numb{cccc},  in  \numb{xxx}  partiti,  octoagies  triginta
dinumerantur.   Duc octoagies  \numb{lxxx}, et  fiunt \numb{vi}  millia  \numb{cccc}.    Tot  in
hujusmodi civitate  domus, secundum  propositionem supra scriptam, construi
[Bede: constitui] possunt.
}{
\probb{Solution.}
This city  measures $8000$  feet around, which is divided into proportions of
one-and-a-half to  one, i.e.  $4800$ and  $3200$.   The length and width of the
houses are  to be  of these  [dimensions].   Thus, take half of each of the
above [measurements],  and from  the larger number there shall remain $2400$,
while from  the the  smaller, $1600$.   Then, divide $1600$ into twenty [parts]
and you will obtain $80$ times $20$.  In a similar fashion, [divide] the larger
number, i.e. $2400$, into $30$ pieces, deriving $80$ times $30$.  Take $80$ times $80$,
making $6400$.   This  many houses  can be  built in  the city, following the
above-written proposal.
}{
%\probb{Решение.}
%В этом городе по окружности $8000$ футов насчитаны, что делится в пропорции полтора [к одному], т.\,е. $4800$ к $3200$.
%Но в первую включена длина, а во вторую ширина домов. 
%Итак, возьми половину от каждого, и от большего числа останется $2400$,
%от меньшего же --- $1600$. Затем раздели $1600$ на $20$ и получишь $80$ раз по $20$ [т.\,е. $1600=80\times 20$],
%аналогичным образом большее число, т.\,е. $2400$, разделенное на $30$, получается $80$ раз по $30$ [т.\,е. $2400=80\times 30$].
%Сделай $80$ раз $80$, получается $6400$. Столько домов можно построить [Bede:  разместить] 
%в этом городе согласно представлениям, написанным выше.
}

\medskip

As Burkholder notes \cite{Burkholder1993},
Alcuin replaces the original round city by a rectangular city with the same circumference $8000$ feet, 
whose sizes are proportional to the house sizes.
The rectangle obtained has sizes $2400\times 1600$ and containes $6400$ houses.
The shortcoming here is that closed curves of equal lenghts can enclose different areas.  

\begin{figure}
\centering
\includegraphics[width=\linewidth]{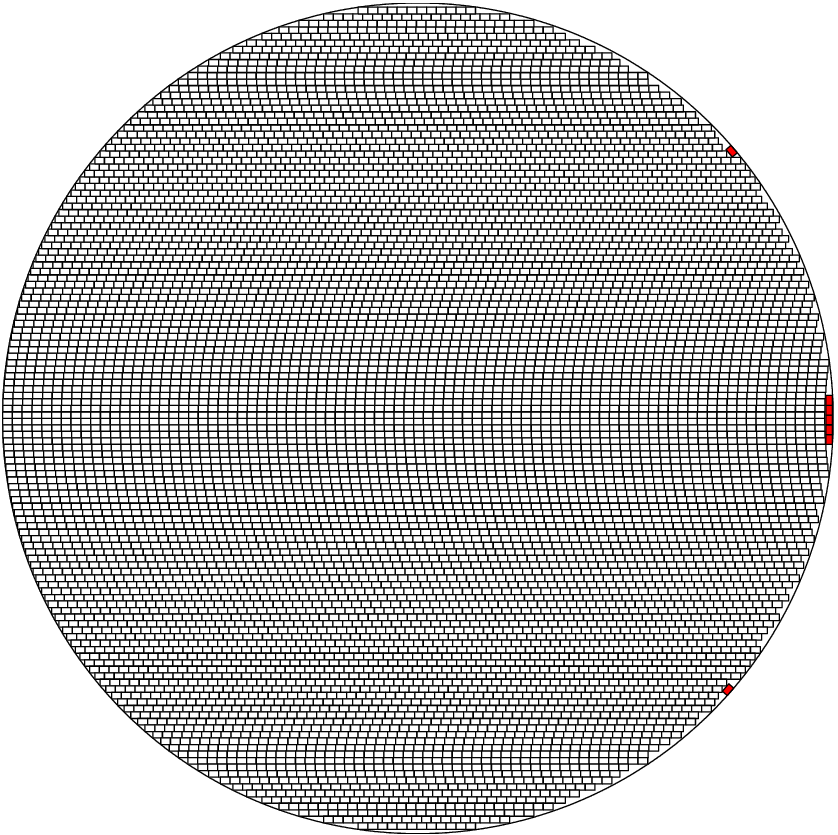} 
\caption{Problem XXVIIII. One can put $8349$ houses in the round city.}
\label{fig_rotunda}
\end{figure}

\begin{figure}
\centering
\includegraphics[width=\linewidth]{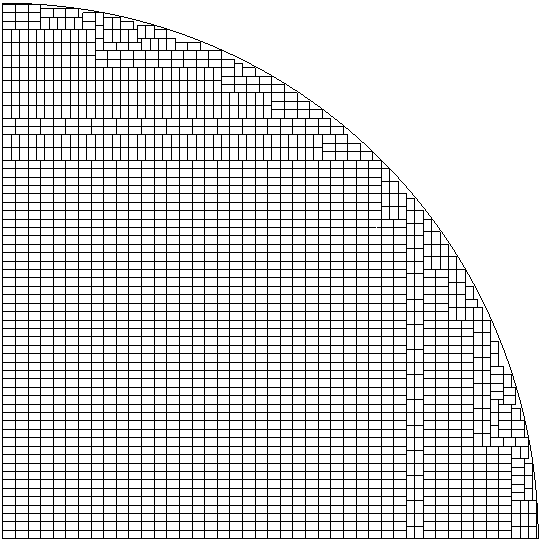} 
\caption{Problem XXVIIII. One can put $8360$ houses in the round city. Here the picture for one quarter is presented \cite{Hemme2010}.}
\label{fig_hemme}
\end{figure}

\medskip

Folkerts \cite{Folkerts1978} gives a second solution:

\medskip

\threel{
Ambitus huius civitatis \numbm{viii} complectitur pedum. Duc ergo quartam de \numbm{viii} partem, fiunt \numbm{ii}. 
Rursusque duc tertiam de \numbm{viii} partem, fiunt \numbm{ii} \numb{dclxvi}.
Duc vero mediam de duobus milibus, fiunt \numbm{i}, atque iterum de duobus milibus \numb{dclxvi} mediam assume partem, fiunt \numbm{i} \numb{cccxxxiii}.
Deinde duc partem tricesimam de \numbm{i} \numb{cccxxxiii}, fiunt \numb{xxxxiiii}, rursusque duc partem vigesimam de \numbm{i}, fiunt \numb{l}. Duc vero quinquagies \numb{xxxxiiii}, fiunt \numbm{ii} \numb{cc}. Deinde duc simul bina milia \numb{cc} quarter, fiunt \numbm{viii} \numb{dccc}. Hoc est summa domorum.}{
This city  measures $8000$  feet around.
Thus, take a thirth quarter of $8000$, making $2000$.
Then  take a $8000$, making $2 666$.
Take a half of $2000$, making $1000$,
and then take a half of $2666$, making $1333$.
Then take a thirtieth part of $1333$, making $44$, then take a twentieth part of $1000$, making $50$.
Take $50$ times $44$, making $2200$.
Take four times $2200$, making $8800$. This is a total number of houses.}{
%Вокруг этого города $8000$ футов.
%Следовательно, произведи четвертую часть от $8000$, получается $2000$.
%Затем произведи третью часть от $8000$, получается $2 666$.
%Произведи же половину от $2000$, получается $1000$,
%а затем такж возьми половину от $2666$, получается $1333$.
%Потом произведи тридцатую часть от $1333$, получается $44$, затем произведи двадцатую часть из $1000$, получается $50$.
%Произведи же $50$ раз по $44$, получается $2200$.
%Потом создай $2200$ четвертый, получается $8800$. Это есть общее количество домов.
}

\smallskip

Here the following formula for the area of a circle is used,
$$
S = \frac{\ell^2}{12}, 
$$
where $\ell$ is the circumference of the circle.
The formula requires $\pi = 3$.
In our case the area could be $5,333,333.333\dots$ sq ft, i.e. $8888.889\dots$ house areas. 
The correct area of the circle is 
$$
S = \frac{\ell^2}{4\pi} = 5,092,958.179\dots \mbox{~sq ft}, 
$$
i.e. $8488.264\dots$ house areas.

Singmaster has fitted in $8307$ houses \cite{HadleySingmaster1992} but he notes that 
``it is probably possible to fit more in''.
I found a solution with $8349$ houses; see Fig.\,\ref{fig_rotunda}.
There are $8342$ houses put ``horizontally'' plus $7$ houses rotated.
H.\,Hemme \cite{Hemme2010} managed to put $8360$ houses!
See Fig.\,\ref{fig_hemme}, where the plan for one quarter of the city is presented.

\end{document}